\documentclass[11pt,reqno,letterpaper]{amsart}
\usepackage{etoolbox}
\usepackage{amsmath,amssymb,mathtools}
\usepackage{microtype}
\usepackage[hidelinks,bookmarksnumbered]{hyperref}
\hypersetup{pdfauthor={Jintian Zhu},pdftitle={On diameter rigidity from two complementary complex submanifolds}}
\numberwithin{equation}{section}
\newtheorem{theorem}{Theorem}[section]
\newtheorem{lemma}[theorem]{Lemma}
\newtheorem{question}[theorem]{Question}
\makeatletter
\patchcmd{\@settitle}{\uppercasenonmath\@title}{}{}{}
\patchcmd{\@settitle}{\bfseries}{\fontsize{14.4}{18}\selectfont\bfseries}{}{}
\makeatother
\newcommand{\CP}{\mathbb{CP}}
\newcommand{\C}{\mathbb C}
\newcommand{\R}{\mathbb R}
\newcommand{\gFS}{g_{\mathrm{FS}}}
\newcommand{\II}{\mathrm{II}}
\DeclareMathOperator{\HSC}{HSC}
\DeclareMathOperator{\Vol}{Vol}

\title[On diameter rigidity from complementary submanifolds]{On diameter rigidity from two complementary complex submanifolds}
\author{Jintian Zhu}
\address[Jintian Zhu]{Institute for Theoretical Sciences, Westlake
University, 600 Dunyu Road, 310030, Hangzhou, Zhejiang,
People's Republic of China}
\email{zhujintian@westlake.edu.cn}


\begin{document}
\begin{abstract}
Let $M$ be a complete connected K\"ahler manifold of complex dimension
$n$ with holomorphic sectional curvature at least $2$. Suppose that $M$
contains two closed, connected complex submanifolds whose complex dimensions sum
to $n-1$ and whose distance is $\pi/\sqrt2$. We prove that $M$ is
biholomorphically isometric to complex projective space with the
Fubini--Study metric of holomorphic sectional curvature $2$. The isometry
carries the two submanifolds onto orthogonal projective linear subspaces.
\end{abstract}
\maketitle

\section{Introduction}

Throughout the paper, we use the convention that the Fubini--Study metric $\gFS$ on $\CP^n$ has
holomorphic sectional curvature $2$ and diameter $$D=\frac{\pi}{\sqrt2}.$$
There has been a broad study on diameter estimates on K\"ahler manifolds.
Tsukamoto~\cite[Theorem~1]{Tsukamoto} proved
that every complete connected K\"ahler manifold $(M^n,g)$ with holomorphic
sectional curvature $\HSC\geq2$ has diameter at most
$D$. For $n\geq2$, Chu, Lee, and the author~\cite[Appendix~B]{CLZ}
showed that
$$\CP^1\times\CP^{n-1}$$ admits a product metric with
$\HSC\geq2$ and diameter $D$.

There are two natural ways for rigidity. One is to enhance the curvature assumption, and the other is to enhance the distance assumption. Along the first direction, under the stronger holomorphic bisectional
curvature bound $\mathrm{BK}\geq1$, Datar and Seshadri~\cite{DS} proved that
if the K\"ahler manifold $(M,g,J)$ has diameter $D$, then it is biholomorphically isometric to
$(\CP^n,\gFS)$. Previous to this work, Liu and Yuan~\cite{LY} partially established this result
under an additional volume assumption.

We focus on the second direction. Note that
 every proper projective linear subspace $\mathbb P(V)$
lies at distance $D$ from its orthogonal complement $\mathbb P(V^\perp)$ in $(\CP^n,\gFS)$. Therefore, it is natural to ask
\begin{question}
Let $(M^n,g,J)$ be a compact connected K\"ahler manifold with
$\HSC\geq2$. Suppose that $M$ contains two nonempty, compact, connected
complex submanifolds $L$ and $N$ satisfying
$$
 \dim_{\C}L+\dim_{\C}N=n-1\,\,\text{and}\,\,d(L,N)=D.
$$
Is $(M,g,J)$ biholomorphically isometric to $(\CP^n,\gFS)$?
\end{question}

Tam and Yu~\cite[Corollary~3.2]{TY} answered this question under the
stronger holomorphic bisectional curvature bound $\mathrm{BK}\geq1$ by
eigenvalue comparison and rigidity in the equality case.
Under $\HSC\geq2$, Chu, Lee, and the author~\cite[Theorem~1.5]{CLZ}
proved the result when one submanifold is a point and the other is a
complex hypersurface. We extend their result to all complementary dimensions.

\begin{theorem}\label{thm:main}
Let $(M^n,g,J)$ be a complete connected K\"ahler manifold with
$\HSC\geq2$. Suppose that $M$ contains two nonempty, closed, connected complex
submanifolds $L$ and $N$ satisfying
$$
 \dim_{\C}L+\dim_{\C}N=n-1\,\,\text{and}\,\,d(L,N)=\frac{\pi}{\sqrt2}.
$$
Then there is a biholomorphic isometry from $(M,g,J)$ to
$(\CP^n,\gFS)$ such that $L$ and $N$ correspond to the projective linear subspaces
$\mathbb P(V)$ and $\mathbb P(V^\perp)$, respectively, for some complex
subspace $V\subset\C^{n+1}$.
\end{theorem}

Fix $q\in N$. By Chu, Lee, and the author~\cite[Proposition~1.6]{CLZ},
each $p\in L$ is antipodal to $q$ on a totally geodesic $\CP^1$.
This gives a surjective endpoint map
$$\varphi_q:\mathcal R_q\to L\,\,\text{with}\,\,\mathcal R_q\subset \mathbb P(\nu_qN),$$
where $\ell\in\mathcal R_q$ if all unit vectors $v\in\ell$ have the same
endpoint $\exp_q(Dv)\in L$.
The map $\varphi_q$ admits local smooth extensions whose differentials are
complex-linear isometries at the points of $\mathcal R_q$.
A volume comparison argument then gives $\mathcal R_q=\mathbb P(\nu_qN)$.
Following the Jacobi field method of Xiong and Yang~\cite{XY}, we recover
the metric in normal coordinates about $N$. Its scalar curvature and
Berger's averaging formula give $\HSC_g\equiv2$, which yields the desired theorem.

\section*{Acknowledgements}
J. Zhu is partially supported by National Key R\&D Program of China
2023YFA1009900, NSFC Grant No.\,12401072, Zhejiang Provincial Natural
Science Foundation of China under Grant No.\,LQKWL26A0101 and the
start-up fund from Westlake University.

\section*{Declaration of AI usage}
The mathematical ideas are due to the author. The author used
ChatGPT (Astra 6.2) to help complete the measure-theoretic argument in
Lemma~\ref{lem:saturation}, check references, and revise the exposition
for clarity. The author takes full responsibility
for the mathematical content.

\section{Proof of the main theorem}

Let $(M^n,g,J)$ denote a complete connected K\"ahler manifold of complex
dimension $n$ with holomorphic sectional curvature at least $2$. Let $L$ and $N$ be two
nonempty, closed, connected complex submanifolds of $M$ such that
$$
 \dim_{\C}L+\dim_{\C}N=n-1\,\,\text{and}\,\,d(L,N)=D.
$$
Along a unit-speed geodesic $\gamma$, we write $T=\gamma'$ and let
$$P_s^t:T_{\gamma(s)}M\to T_{\gamma(t)}M$$ denote parallel transport from $\gamma(s)$ to $\gamma(t)$.

\begin{lemma}\label{lem:radial}
Every unit-speed minimizing geodesic $\gamma$ from $q\in N$ to $p\in L$
has length $D$ and meets both submanifolds orthogonally. Along $\gamma$,
there is an orthogonal decomposition
$$
 T_{\gamma(t)}M=\C T\oplus P_0^t(T_qN)\oplus P_D^t(T_pL).
$$
The radial curvature operator satisfies
\begin{equation}\label{eq:operator}
 R(JT,T)T=2JT\,\,\text{and}\,\,
 R(E,T)T=\tfrac12 E\,\,\text{for any }E\perp\C T.
\end{equation}
The shape operator of $N$ in the direction $T(0)$ and that of $L$ in
the direction $-T(D)$ both vanish.
\end{lemma}

\begin{proof}
By Tsukamoto~\cite[Theorem~1]{Tsukamoto}, $M$ is compact and
$\operatorname{diam}(M,g)\leq D$. Hence
\begin{equation}\label{eq:endpoint-distance}
 D=d(N,L)\leq d(q,p)\leq\operatorname{diam}(M,g)\leq D.
\end{equation}
In particular, $\gamma$ minimizes the distance between $N$ and $L$. The first
variation formula gives $T(0)\perp T_qN$ and $T(D)\perp T_pL$. Since $T$ and $JT$ are parallel
along $\gamma$, $P_0^t(T_qN)$ and $P_D^t(T_pL)$ are orthogonal to $\C T$.

By \cite[Remark~2.2]{CLZ}, $\sin(\sqrt2 t)JT$ is a Jacobi field.
Since $JT$ is parallel, the Jacobi equation combined with the continuity gives
$$R(JT,T)T=2JT$$ on $[0,D]$.
In terms of the sectional curvature $K$, the bisectional curvature
bound in \cite[Lemma~2.1]{CLZ} reads
$$
 K(T,E)+K(T,JE)\geq1\,\,\text{for every unit }E\perp\C T.
$$

For a unit vector $X\in T_qN$, set
$$E(t)=P_0^tX\,\,\text{and}\,\,f(t)=\cos(t/\sqrt2).$$
Since $f(D)=0$, both $fE$ and $fJE$ are admissible for variations from
$N$ to $L$. Let $I$ denote the index form for these variations.
Since the second fundamental form
of the complex submanifold $N$ satisfies $$\II^N(JX,JX)=-\II^N(X,X),$$
the second variation formula and the curvature bound give
$$
 \begin{aligned}
 0&\leq I(fE,fE)+I(fJE,fJE)\\
  &=\int_0^D\bigl(2(f')^2-f^2(K(T,E)+K(T,JE))\bigr)\,dt\\
  &\leq\int_0^D\bigl(2(f')^2-f^2\bigr)\,dt
   =0 .
 \end{aligned}
$$
This implies $I(fE,fE)=0$, and for every
admissible $V$ we have
$$I(fE,V)=\frac{1}{2}\left.\frac{d}{dt}\right|_{t=0}I(fE+tV,fE+tV)=0.$$ Integrating by parts and using $f''=-f/2$ gives
$$
 \begin{aligned}
 0={}&-\int_0^D f\big\langle R(E,T)T-\tfrac12 E,V\big\rangle\,dt\\
     &-\tfrac1{\sqrt2}\langle P_0^D X,V(D)\rangle
      -\langle\II^N(X,V(0)),T(0)\rangle .
 \end{aligned}
$$
First take $V$ supported in $(0,D)$. Since $f>0$, we have
$$R(E,T)T=E/2$$ on $(0,D)$, which is also extended on $[0,D]$ by continuity.
Taking $V(D)=0$ and $V(0)\in T_{q}N$ arbitrary gives $$\langle\II^N(X,Z),T(0)\rangle=0$$ for every
$Z\in T_qN$. Taking $V(0)=0$ and $V(D)\in T_pL$ arbitrary gives $$P_0^D X\perp T_pL.$$

By linearity the identities above hold for every $X\in T_qN$.
Thus the shape operator of $N$ in the direction
$T(0)$ vanishes, and $P_0^t(T_qN)$ is orthogonal to $P_D^t(T_pL)$
for every $t$. Reversing $\gamma$ gives the same curvature identity
on $P_D^t(T_pL)$ and the vanishing of the shape operator of $L$ in
the direction $-T(D)$. From the dimension counting argument, we obtain the desired decomposition for $T_{\gamma(t)}M$ and also the curvature identity \eqref{eq:operator}.
\end{proof}

Fix $q\in N$ and let $\nu_qN$ be the normal space to $N$ at $q$.
We equip the projective normal space $\mathbb P(\nu_qN)$ with the Fubini--Study
metric. Let $\mathcal R_q$ denote
the set of complex lines $\ell\in \mathbb P(\nu_qN)$ for which
$\exp_q(Dv)$ is the same point of $L$ for every unit vector $v\in\ell$.
Denote this common endpoint by $\varphi_q(\ell)$.

\begin{lemma}\label{lem:endpoint}
The set $\mathcal R_q$ is compact and $\varphi_q:\mathcal R_q\to L$ is
surjective. Near each $\ell\in\mathcal R_q$, the map $\varphi_q$ extends
to a diffeomorphism from a neighborhood $U$ of $\ell$ in $\mathbb P(\nu_qN)$ onto an
open subset of $L$. Its differential is a complex-linear isometry at
every point of $U\cap\mathcal R_q$.
\end{lemma}

\begin{proof}
      Let $S(\nu_qN)$ denote the unit sphere in $\nu_qN$.
Define 
$$h(v)=\exp_q(Dv)$$ for any $v\in S(\nu_qN)$, and set
$$
 \widetilde{\mathcal R}_q
 =\bigl\{v\in S(\nu_qN)\mid h(e^{\theta J}v)=h(v)\in L
 \,\,\text{for every }\theta\in\R\bigr\}.
$$
Since $h$ is continuous and $L$ is closed, $\widetilde{\mathcal R}_q$
is closed in $S(\nu_qN)$ and hence compact. Its image under the
Hopf map $v\mapsto\C v$ is $\mathcal R_q$, so $\mathcal R_q$ is compact.

Let $p\in L$. By \eqref{eq:endpoint-distance}, we have $d(q,p)=D$. By
\cite[Proposition~1.6]{CLZ}, $q$ and $p$ are antipodes on a totally
geodesic complex curve isometric to a round sphere. Each meridian
from $q$ to $p$ has length $D$ and so attains the distance between $N$ and $L$. This implies that its tangent vector at $q$
is normal to $N$. Therefore, the complex tangent line of the complex curve at $q$ belongs
to $\mathcal R_q$ and maps to $p$ by $\varphi_q$. This establishes the surjectivity of $\varphi_q$.

To construct the extension, choose a local unit section $\sigma$ of
the Hopf map in a small neighborhood of $\ell\in\mathcal R_q$ and project $h\circ\sigma$
normally onto $L$. This gives a smooth map $\Phi$ near $\ell$ which
agrees with $\varphi_q$ on $\mathcal R_q$. Let us show that the differential map of this extended map is a complex-linear isometry at
 any $z\in\mathcal R_q$ in this neighborhood, then the lemma follows after passing to a smaller neighborhood by the inverse function theorem.
 
 Put $v=\sigma(z)$
and $p=h(v)$. Define $$\gamma(t)=\exp_q(tv), \,0\leq t\leq D.$$
 By Lemma~\ref{lem:radial}, parallel
transport along $\gamma$ identifies the horizontal space
$$\mathcal H=\nu_qN\cap(\C v)^\perp$$ with $T_pL$.
For any horizontal vector $W\in \mathcal H$, the Jacobi field $V(t)$ with initial data
$(V(0),V'(0))=(0,W)$ is $V(t)=\sqrt2\sin(t/\sqrt2)P_0^tW$ by \eqref{eq:operator}, and so we have
$$dh_v(W)=\sqrt2 P_0^DW.$$ 
We have $dh_v(Jv)=0$,
since $h$ is constant on the Hopf circle through $v$.

Let $W\in \mathcal H$ be the horizontal lift of $\zeta\in T_z\mathbb P(\nu_qN)$. Then
$d\sigma_z(\zeta)-W$ is a multiple of $Jv$. Since the normal projection
fixes $L$ and $P_0^DW\in T_pL$, we obtain
$$
 d\Phi_z(\zeta)=\sqrt2 P_0^DW .
$$
Note that the horizontal lift and parallel transport commute with $J$, and we have
$$|\zeta|_{\gFS}^2=2|W|^2.$$ 
We conclude that $$d\Phi_z: T_z\mathbb P(\nu_qN)\to T_pL$$ is a complex-linear
isometry as desired.
\end{proof}

Let $g_N$ and $g_L$ denote the metrics induced by $g$ on $N$ and $L$,
respectively. We next identify $(L,g_L)$ with the projective normal
space at each point of $N$.

\begin{lemma}\label{lem:saturation}
For every $q\in N$, we have $\mathcal R_q=\mathbb P(\nu_qN)$, and the
endpoint map
$$
 \varphi_q:(\mathbb P(\nu_qN),\gFS)\longrightarrow(L,g_L)
$$
is a biholomorphic isometry. Similarly, for every $p\in L$, the endpoint map
$$
 \varphi_p:(\mathbb P(\nu_pL),\gFS)\longrightarrow(N,g_N)
$$
is a biholomorphic isometry. Both $N$ and $L$ are totally geodesic in $M$.
\end{lemma}

\begin{proof}
Fix $q\in N$. By Lemma~\ref{lem:endpoint} and compactness, we can choose finitely
many open balls $B_i$ in $\mathbb P(\nu_qN)$ to cover $\mathcal R_q$,
with each $\overline{B_i}$ contained in a neighborhood on which
$\varphi_q$ extends to a diffeomorphism $\Phi_i$ onto an open subset
of $L$. Set $$K_i=\mathcal R_q\cap\overline{B_i}.$$
Then $\Phi_i|_{K_i}=\varphi_q|_{K_i}$, and $(d\Phi_i)_\ell$ is a
complex-linear isometry for every $\ell\in K_i$.

By surjectivity of $\varphi_q$, the compact sets $\varphi_q(K_i)$ cover
$L$. Write $$O_i=\operatorname{int}_L\varphi_q(K_i).$$ Since
$$
 L\setminus\bigcup_i O_i\;\subset\;\bigcup_i\partial_L\varphi_q(K_i)
$$
and each boundary is closed and nowhere dense in $L$, $\bigcup_i O_i$
is dense in $L$. For $p\in O_i$, we have $\Phi_i^{-1}(p)\in K_i$,
so $(d\Phi_i^{-1})_p$ is a complex-linear isometry. Hence
$$\Phi_i^{-1}:O_i\to\Phi_i^{-1}(O_i)$$ is a biholomorphic isometry onto
an open subset of $(\mathbb P(\nu_qN),\gFS)$. This yields
$\HSC_{g_L}(X)=2$ for every unit vector $X\in T_pL$ with
$p\in\bigcup_i O_i$. By continuity on the unit tangent bundle of $L$,
the equality holds for every $p\in L$ and every unit $X\in T_pL$.

Since $L$ is closed in $M$, it is compact. In positive dimension,
Tsukamoto~\cite[Theorem~3]{Tsukamoto} gives simple connectedness of $L$,
and the classification of complete simply connected complex space forms
\cite[Chapter~IX, \S7]{KN} gives an isometry from $L$ to $\mathbb P(\nu_qN)$.
If $\dim_{\C}L=0$, both spaces are singletons. In particular,
\begin{equation}\label{eq:projective-normal-volume}
 \Vol(L)=\Vol(\mathbb P(\nu_qN)).
\end{equation}

The Borel sets $E_i=K_i\setminus\bigcup_{j<i}K_j$ partition $\mathcal R_q$.
On $E_i$, we have $\varphi_q=\Phi_i$ and the Jacobian of $\Phi_i$ is one.
Since the sets $\Phi_i(E_i)$ cover $L$, change of variables gives
$$
 \Vol(L)\leq\sum_i\Vol(\Phi_i(E_i))
 =\sum_i\Vol(E_i)=\Vol(\mathcal R_q)
 \leq\Vol(\mathbb P(\nu_qN)).
$$
By \eqref{eq:projective-normal-volume}, the open set
$\mathbb P(\nu_qN)\setminus\mathcal R_q$ has zero volume and is therefore
empty. By Lemma~\ref{lem:endpoint}, $\varphi_q:\mathbb P(\nu_qN)\to L$
is now a surjective local holomorphic isometry. Its covering degree is
one by \eqref{eq:projective-normal-volume}, so $\varphi_q$ is a
biholomorphic isometry.

Every unit normal vector to $N$ is therefore the initial velocity of a
minimizing geodesic to $L$. By Lemma~\ref{lem:radial}, the shape operator of $N$ vanishes in
every normal direction, so $N$ is totally geodesic. Interchanging $N$ and
$L$ gives the endpoint isometry onto $N$ and shows that $L$ is totally geodesic.
\end{proof}

On the unit normal bundle $S(\nu N)$, define
$$
 \pi_N(q,v)=q\,\,\text{and}\,\,\pi_L(q,v)=\exp_q(Dv).
$$
For the velocity $\xi$ of a curve $(q(s),v(s))$ in $S(\nu N)$, set
$$\eta(\xi)=\langle(\nabla^\perp_s v)(0),Jv\rangle.$$

\begin{lemma}\label{lem:circle}
The map $$(\pi_N,\pi_L):S(\nu N)\to N\times L$$ is a principal circle
bundle with action $v\mapsto e^{\theta J}v$ and connection form $\eta$.
\end{lemma}

\begin{proof}
Define the projection
$$
 \rho:S(\nu N)\to\mathbb P(\nu N),\,\, (q,v)\mapsto(q,\C v).
$$
The map $\rho$ is the unit circle bundle of the tautological line
bundle, with circle action $v\mapsto e^{\theta J}v$.
By Lemma~\ref{lem:saturation}, we have $(\pi_N,\pi_L)=\Psi\circ\rho$, where
$$
 \Psi:\mathbb P(\nu N)\to N\times L,\,\,
 (q,\ell)\mapsto(q,\varphi_q(\ell)).
$$
Near each $(q,\ell)\in\mathbb P(\nu N)$, choose a smooth local section
$\sigma$ of the Hopf projection $\rho$. On the domain of $\sigma$, we
have $\Psi=(\pi_N,\pi_L)\circ\sigma$. Hence $\Psi$ is smooth.
Each $\varphi_q$ is a diffeomorphism by Lemma~\ref{lem:saturation},
so $\Psi$ is bijective. In local bundle coordinates we have
$$
 d\Psi_{(q,\ell)}=
 \begin{pmatrix}
  \mathrm{Id}_{T_qN} & 0\\
  * & (d\varphi_q)_\ell
 \end{pmatrix}.
$$
Both diagonal blocks are invertible, so $\Psi$ is a diffeomorphism
and identifies $(\pi_N,\pi_L)$
with the principal circle bundle $\rho:S(\nu N)\to\mathbb P(\nu N)$.

The normal connection preserves $J$, so for each fixed $\theta$ we have
$\nabla_s^\perp(e^{\theta J}v)=e^{\theta J}\nabla_s^\perp v$.
Since $e^{\theta J}$ preserves the metric, $\eta$ is invariant under
the circle action by definition. On the other hand, the orbit $s\mapsto e^{sJ}v$ has velocity $Jv$
at $s=0$, and $\eta(Jv)=\langle Jv,Jv\rangle=1$.
Hence $\eta$ is a connection form.
\end{proof}

\begin{lemma}\label{lem:normal-metric}
The normal exponential map
$$
 \mathcal E:S(\nu N)\times(0,D)\to M\setminus(N\cup L),\,\,
 (q,v,t)\mapsto\exp_q(tv),
$$
is a diffeomorphism, and
\begin{equation}\label{eq:tube-metric}
 \mathcal E^*g=dt^2+\cos^2(t/\sqrt2)\,\pi_N^*g_N+\sin^2(t/\sqrt2)\,\pi_L^*g_L
 +\tfrac12\sin^2(\sqrt2 t)\,\eta^2.
\end{equation}
In this formula $\pi_N$ and $\pi_L$ are composed with the projection
$(u,t)\mapsto u$, and $\eta^2=\eta\otimes\eta$ is pulled back by the
same projection.
\end{lemma}

\begin{proof}
Every unit-speed geodesic starting normally from $N$ reaches $L$ at
time $D$ by Lemma~\ref{lem:saturation}. Since $d(N,L)=D$, each such
geodesic minimizes from $N$ to $L$ on $[0,D]$. Given
$x\in M\setminus(N\cup L)$, a shortest geodesic from $N$ to $x$
starts normally and has length $d(N,x)\leq D$ by
\eqref{eq:endpoint-distance}. Equality would give $x\in L$.
Thus $0<d(N,x)<D$ and $\mathcal E$ is surjective.

If two unit normal geodesics from $N$ pass through $x$ before time $D$,
both reach $x$ at time $d(N,x)$. Following the first from $N$ to $x$
and the second from $x$ to $L$ gives a curve of length $D$. The curve
minimizes from $N$ to $L$, so the two velocities agree at $x$.
Uniqueness for the geodesic equation then gives injectivity of $\mathcal E$.

Thus $\mathcal E$ is a smooth bijection between manifolds of the same
dimension. To prove that $\mathcal E$ is a diffeomorphism, it suffices
to show that $\mathcal E^*g$ is positive definite.
For $\xi\in T_{(q,v)}S(\nu N)$, choose a curve $(q(s),v(s))$ with
initial velocity $\xi$. The Jacobi field along $\gamma(t)=\exp_q(tv)$ given by
$$
 V_\xi(t)=\left.\frac{\partial}{\partial s}\right|_{s=0}
             \exp_{q(s)}(tv(s))
$$
has endpoint values $V_\xi(0)=d\pi_N(\xi)$ and $V_\xi(D)=d\pi_L(\xi)$,
both orthogonal to $\C T$ by Lemma~\ref{lem:radial}. Thus
$\langle V_\xi,T\rangle$, being affine with zero endpoint values, vanishes.
In a parallel frame on $(\C T)^\perp$, each coefficient $f$ satisfies
$f''+f/2=0$ by \eqref{eq:operator}, and so we have
$$f(t)=f(0)\cos(t/\sqrt2)+f(D)\sin(t/\sqrt2).$$
On the other hand, the coefficient $h=\langle V_\xi,JT\rangle$ satisfies $h''+2h=0$,
$h(0)=0$ and $h'(0)=\langle(\nabla_s v)(0),Jv\rangle=\eta(\xi)$.
Therefore, we have
\begin{equation}\label{eq:polar}
 \begin{aligned}
 V_\xi(t)={}&\cos(t/\sqrt2)\,P_0^t\bigl(d\pi_N(\xi)\bigr)
       +\sin(t/\sqrt2)\,P_D^t\bigl(d\pi_L(\xi)\bigr)\\
       &+\frac{\sin(\sqrt2 t)}{\sqrt2}\,\eta(\xi)JT.
 \end{aligned}
\end{equation}
By Lemma~\ref{lem:radial}, the three summands in \eqref{eq:polar} are
mutually orthogonal and lie in $T^\perp$. Since $d\mathcal E(\partial_t)=T$
and $d\mathcal E(\xi,0)=V_\xi(t)$, the radial part of $\mathcal E^*g$
is $dt^2$ and the mixed terms vanish. Taking the squared norm of
$V_\xi(t)$ in \eqref{eq:polar} gives the remaining terms in
\eqref{eq:tube-metric}.

For $0<t<D$, all coefficients in \eqref{eq:tube-metric} are positive.
If $V_\xi(t)=0$, then $d\pi_N(\xi)=d\pi_L(\xi)=0$ and $\eta(\xi)=0$.
By Lemma~\ref{lem:circle}, $\ker d(\pi_N,\pi_L)=\R Jv$ and $\eta(Jv)=1$,
so $\xi=0$. Thus $\mathcal E^*g$ is positive definite, and the smooth
bijection $\mathcal E$ is a diffeomorphism.
\end{proof}

\begin{proof}[Proof of Theorem~\ref{thm:main}]
We first compute the real scalar curvature of $g$ on $M\setminus(N\cup L)$.
Define the circle bundle projection
$$
 \pi:S(\nu N)\times(0,D)\longrightarrow N\times L\times(0,D),\,\,
 (u,t)\longmapsto(\pi_N(u),\pi_L(u),t).
$$
By \eqref{eq:tube-metric}, we have $\mathcal E^*g=\pi^*h+r^2\eta^2$, where
$r(t)=\sin(\sqrt2 t)/\sqrt2$ and
$$
 h=dt^2+\cos^2(t/\sqrt2)\,\pi_N^*g_N+\sin^2(t/\sqrt2)\,\pi_L^*g_L.
$$
On the base, $\pi_N(q,p,t)=q$ and $\pi_L(q,p,t)=p$.
Let $U$ be the unit vector in the direction of the circle action, and
let $\Omega$ be the curvature form on the base, so that $\pi^*\Omega=d\eta$.
With $\Delta_h r=\operatorname{div}_h(\nabla_h r)$, the submersion formulas
of Besse~\cite[Proposition~9.36]{Besse}, in the circle bundle form
of \cite[(66), (68)]{BNS}, give
$$
 \begin{aligned}
 \operatorname{Scal}_{\mathcal E^*g}
 &=\pi^*\!\left(\operatorname{Scal}_h-2\frac{\Delta_h r}{r}
               -\frac{r^2}{2}|\Omega|_h^2\right),\\
 \operatorname{Ric}_{\mathcal E^*g}(U,U)
 &=\pi^*\!\left(-\frac{\Delta_h r}{r}
               +\frac{r^2}{2}|\Omega|_h^2\right).
 \end{aligned}
$$
Adding the two identities gives
$$
 \operatorname{Scal}_{\mathcal E^*g}
 +\operatorname{Ric}_{\mathcal E^*g}(U,U)
=\pi^*\!\left(\operatorname{Scal}_h-3\,\frac{\Delta_h r}{r}\right).
$$
The circle orbit $s\mapsto(q,e^{sJ}v,t)$ has velocity $(Jv,0)$ at $s=0$, and \eqref{eq:polar} gives $d\mathcal E(Jv,0)=r(t)JT$.
Since $(Jv,0)$ is vertical and $\eta(Jv)=1$, \eqref{eq:tube-metric}
gives $|(Jv,0)|_{\mathcal E^*g}=r(t)$. Thus
$$
 U=\frac{(Jv,0)}{r(t)}\,\,\text{and}\,\, d\mathcal E(U)=JT.
$$
The K\"ahler symmetry of the Ricci tensor and \eqref{eq:operator} give
$$
 \operatorname{Ric}_{\mathcal E^*g}(U,U)
 =\operatorname{Ric}_g(JT,JT)=\operatorname{Ric}_g(T,T)
 =2+\tfrac12(2n-2)=n+1.
$$

Put $k=\dim_{\C}N$, $\ell=\dim_{\C}L$ and $s=t/\sqrt2$, so that
$k+\ell=n-1$. By Lemma~\ref{lem:saturation},
$\operatorname{Scal}_{g_N}=2k(k+1)$ and
$\operatorname{Scal}_{g_L}=2\ell(\ell+1)$. Substituting the warping
functions $\cos s$ and $\sin s$ and the real fiber dimensions $2k$ and
$2\ell$ into the scalar curvature formula of Dobarro and
\"Unal~\cite[Proposition~2.6(1)]{DU} gives
$$
 \operatorname{Scal}_h=(2n-1)(n-1)
       +\frac{3k}{\cos^2 s}+\frac{3\ell}{\sin^2 s}.
$$
For the radial function $r$, \cite[Proposition~2.3(3)]{DU} gives
$$
 \frac{\Delta_h r}{r}
 =\frac{\bigl(\cos^{2k}s\,\sin^{2\ell}s\,r'\bigr)'}
        {\cos^{2k}s\,\sin^{2\ell}s\,r}
 =-2n+\frac{k}{\cos^2 s}+\frac{\ell}{\sin^2 s}.
$$
It follows that $\operatorname{Scal}_h-3\Delta_h r/r=(2n+1)(n+1)$.
Subtracting $\operatorname{Ric}_{\mathcal E^*g}(U,U)=n+1$ gives
$\operatorname{Scal}_g=2n(n+1)$ on $M\setminus(N\cup L)$ and hence
on $M$ by continuity.

By Berger's formula (see Kobayashi~\cite[(7.6.44)]{Kobayashi}), the average of $\HSC_g$
over the unit sphere in $T_pM$ is $\operatorname{Scal}_g(p)/(n(n+1))=2$
for every $p\in M$.
Since $\HSC_g\geq 2$, we have $\HSC_g\equiv 2$.
The manifold $M$ is simply connected by \cite[Theorem~3]{Tsukamoto}.
The classification of complete simply
connected complex space forms \cite[Chapter~IX, \S7]{KN} now gives a
biholomorphic isometry $\Phi:(M,g)\to(\CP^n,\gFS)$.

Let $V\subset\C^{n+1}$ be the complex span of the lines representing
the points of $\Phi(L)$. By \eqref{eq:endpoint-distance}, every point
of $\Phi(N)$ lies at distance $D$ from every point of $\Phi(L)$, since two
points of $(\CP^n,\gFS)$ are at distance $D$ precisely when their
representing lines are orthogonal. Thus we obtain
$$
 \Phi(L)\subset\mathbb P(V)\,\,\text{and}\,\, \Phi(N)\subset\mathbb P(V^\perp).
$$
The two projective spaces have complex dimensions summing to
$n-1=\dim_{\C}L+\dim_{\C}N$, so both inclusions have codimension zero.
Since $\Phi(L)$ and $\Phi(N)$ are closed and the projective spaces are
connected, both inclusions are equalities.
\end{proof}

\end{document}